\documentclass[12pt]{amsart}
\usepackage[margin=1.0in]{geometry}               

\usepackage{enumerate}
\usepackage{amssymb}
\usepackage{amsmath}
\usepackage{amscd}
\usepackage{amsthm}
\usepackage{amsfonts}
\usepackage{graphicx}
\usepackage[all]{xy}
\usepackage{verbatim}
\usepackage{hyperref}
\usepackage{epstopdf}
\usepackage[utf8]{inputenc}
\usepackage{bbm}
\usepackage{pdfsync}

\theoremstyle{definition}

\theoremstyle{definition}

\theoremstyle{definition}

\theoremstyle{definition}
\theoremstyle{definition}
\theoremstyle{definition}

\newtheorem{theo}{Theorem}[section]
\newtheorem{prop}{Proposition}[section]
\newtheorem{lemma}{Lemma}[section]
\newtheorem{cor}{Corollary}[section]

\newcommand{\al}{\alpha}

\newcommand{\cA}{\mathcal{A}}
\newcommand{\cB}{\mathcal{B}}
\newcommand{\cC}{\mathcal{C}}

\newcommand{\cP}{\mathcal{P}}

\newcommand{\cT}{\mathcal{T}}

\newcommand{\bG}{\mathbb{G}}

\newcommand{\bR}{\mathbb{R}}
\newcommand{\bZ}{\mathbb{Z}}
\newcommand{\bQ}{\mathbb{Q}}

\newcommand{\bN}{\mathbb{N}}
\newcommand{\bH}{\mathbb{H}}

\newcommand{\gog}{\mathfrak{g}}
\newcommand{\gof}{\mathfrak{f}}
\newcommand{\goh}{\mathfrak{h}}

\newcommand{\gom}{\mathfrak{m}}

\newcommand{\GL}{\operatorname{GL}}

\newcommand{\Mat}{\operatorname{Mat}}

\renewcommand{\hom}{\Gamma\backslash G}

\newcommand\set[1]{\left\{#1\right\}} 

\newcommand\on[1]{\operatorname{#1}} 
\newcommand\diag[1]{\operatorname{diag}\left(#1\right)}

\newcommand{\acts}{\hspace{-1pt}\mbox{\raisebox{1.3pt}{\text{\huge{.}}}}\hspace{-1pt}}

\newcommand{\onto}{\xymatrix{\ar@{>>}[r]&}}
\newcommand{\da}[4]{\xymatrix{#1 \ar@<.5ex>[r]^{#2} \ar@<-.5ex>[r]_{#3} & #4}}

\newif\ifdraft\drafttrue

\theoremstyle{theorem}
\theoremstyle{theorem}

\newcommand{\homq}{\Gamma\backslash G}
\newcommand{\SLQ}[1][m]{\operatorname{SL}_{#1}(\bQ_p)}
\newcommand{\SLR}[1][m]{\operatorname{SL}_{#1}(\bR)}
\newcommand{\SLZ}[1][m]{\operatorname{SL}_{#1}(\bZ_p)}
\newcommand{\pk}{\cP^X_k}
\renewcommand{\al}[1][-l]{a^{#1}\acts}

\newcommand{\sh}{{G^{-}_a}}
\newcommand{\uh}{{G^{+}_a}}
\newcommand{\sha}{{\gog^{-}}}
\newcommand{\uha}{{\gog^{+}}}
\newcommand{\za}{{\gog^{0}}}

\newcommand{\slq}[1][d]{\mathfrak{sl}_{#1}(\mathbb{Q}_{p})}

\newcommand{\Pk}{{K^{P}_{k}}}

\def\bprop{\begin{prop}}
\def\eprop{\end{prop}}
\def\blem{\begin{lema}}
\def\elem{\end{lema}}
\def\bth{\begin{theo}}
\def\eth{\end{theo}}
\def\bcor{\begin{cor}}
\def\ecor{\end{cor}}
\def\proc#1{\medbreak\noindent{\it #1}\hspace{1ex}\ignorespaces}
\def\ep{\noindent{\hfill $\Box$}}

\def\ben{\begin{enumerate}}
\def\een{\end{enumerate}}

\def\acks{\subsubsection*{Acknowledgements.}}

\title[Effectivity of Uniqueness of the Maximal Entropy Measure]{Effectivity of Uniqueness of the Maximal Entropy Measure on $p$-adic homogeneous spaces}
\author[Rene Rühr]{Rene Rühr\\ \\Eidgen\"ossische Technische Hochschule Z\"urich}

\begin{document}
\maketitle

\begin{abstract}
We consider the dynamical system given by an $\on{Ad}$-diagonalizable element $a$ of the $\bQ_p$-points $G$ of a unimodular linear algebraic group acting by translation on a finite volume quotient $X$. Assuming that this action is exponentially mixing (e.g.\ if $G$ is simple) we give an effective version (in terms of $K$-finite vectors of the regular representation) of the following statement: If $\mu$ is an $a$-invariant probability measure with measure-theoretical entropy close to the topological entropy of $a$, then $\mu$ is close to the unique $G$-invariant probability measure of $X$.
\end{abstract}

\section{Statement of the Result}
Let $\bG$ be an algebraic subgroup of $\on{SL}_{m}$ defined over $\bQ$ such that $G=\bG(\bQ_{p})$ is a $d$-dimensional unimodular group. Suppose further that $\Gamma$ is a lattice of $G$ and let us consider the (necessarily compact) quotient $X=\hom$ with unique right-$G$-invariant probability measure $m_X$ which we call Haar measure on $X$. We normalize the Haar measure $m_G$ such that it is compatible with $m_X$.
Let $a\in G$ be $\on{Ad}$-diagonalizable over $\bQ_{p}$ such that at least one eigenvalue is not of $p$-adic norm one and consider the corresponding map on $X$, $x\mapsto a\acts x=xa^{-1}$. Denote by $h_{\mu}(a)$ the measure-theoretic entropy of an $a$-invariant Borel probability measure $\mu$ on $X$. $h_{m_X}(a)$ equals the topological entropy $h_{\text{top}}(a)$ of $a$.

Assume that the group $\langle G^{+}_{a},G^{-}_{a}\rangle$ generated by the horospherical groups $G^{+}_{a}$ and $G^{-}_{a}$ acts uniquely ergodically on $X$, that is, the Haar measure $m_{X}$ is the unique measure invariant under $\langle G^{+}_{a},G^{-}_{a}\rangle$. Then $m_{X}$ is also the unique measure of maximal entropy with respect to translation by $a$. This goes back to \cite{adler} for toral automorphisms and see e.g.\ \cite{margulis}, \cite{pisa} for the present case.
Call $f:X\to\bR$ smooth if it is locally constant and attach the integer
\begin{align*}
l_{f}=\min\{l\in\bN\;|\; f(xg)=f(x)\; &\text{ for all } x\in X \text{ and for all } g\in G \text{ s.t. } g\equiv e \;(\on{mod}{p^{l}})\}
\end{align*}
to it. Equivalently, $f$ is a $K=\SLZ\cap G$-finite vector of the right-regular representation $L^{2}(X,m_{X})$ of $G$ ramified of level $l_{f}$. 
We want to consider an effective analogue of unique ergodicity of $\langle G^{+}_{a},G^{-}_{a}\rangle$, namely the following mixing assumption which we impose: The action by $a$ is exponentially fast mixing, that is, there exist strictly positive constants $c,\alpha,\delta$ such that for any two smooth vectors $f,g\in L^{2}_{0}(X)$ of vanishing integral and ramified of level $l$ we have 
\begin{equation}
\label{assumption}
|\langle f\circ a^n,g\rangle\leq c p^{l\alpha}\|f\|_{L^{2}(X)}\|g\|_{L^{2}(X)}\|a\|^{-\delta n} \quad\text{ for all }n\in\bZ.
\end{equation}
We note that this is always the case if $G$ is simple.
\begin{theo}
\label{theorem}
Let $a\in G$ act on $X$ with the assumption made above. Then there exists a constant $\kappa=\kappa(a,\alpha,\delta,c,X)$ such that for any $a$-invariant Borel probability measure $\mu$ and smooth function $f$, it holds that
$$\left|\int fdm_{X}-\int f d\mu\right| \leq \kappa p^{(2\alpha+\frac{d}{2})l_{f}}\|f\|_{L^{2}(X,m_{X})}(h_{m_{X}}(a)-h_\mu(a))^{\frac{1}{2}}.$$
\end{theo}
The constant $\kappa$ is explicitly calculated to be $$\kappa=\sqrt{2}cp^{2\alpha}m_G(K^G_{2})^{-\tfrac12}(1-\|a\|^{\delta})^{-1}\exp{\left((3\alpha+d)h_{m_{X}}(a)\right)}$$ where $K^{G}_{2}$ denotes the ball of radius $p^{-2}$ in $G$ and $\alpha$, $\delta$ and $c$ are the constants from line (\ref{assumption}).

Analogous results have been obtained for toral automorphisms and for hyperbolic maps on Riemann manifolds in \cite{Polo} and \cite{shirali} respectively. We comment in the next chapter on the adoption for real quotients.

\section{Outline of the Proof}
\label{methodofproof}
Denote $\sh$ and $\uh$ the stable and unstable horospherical subgroup of $G$ with respect to $a$:
$$\sh=\{g\in G: a^{n}ga^{-n}\to e \text{ as } n\to\infty\}$$
$$\uh=\{g\in G: a^{n}ga^{-n}\to e \text{ as } n\to-\infty\}.$$ 
The assumption on $a$ having an eigenvalue not of absolute value one implies that there is at least one contracting and one expanding direction in the Lie algebra so that $\sh$ is non-trivial. Let $$\operatorname{mod}(a,\sh)=|\det{\operatorname{Ad}_{a}|_{\operatorname{Lie}(\sh)}}|$$ denote the modular character corresponding to the conjugation action of $a$ restricted to $\sh$. We recall that in this notation $$h_{\text{top}}(a)=\log{\operatorname{mod}(a,\sh)}$$ whose calculation we provide in Proposition \ref{bowenentropy}. One can formulate the following variational principle
\begin{theo}[Theorem 7.6, Theorem 7.9 of \cite{pisa}]
\label{pisathm}
For any $a\in G$ and any $a$-invariant probability measure $\mu$ on $X$ the entropy of $\mu$ is bounded by
$$h_{\mu}(a)\leq \log{\operatorname{mod}(a,\sh)}$$
and equality holds if and only if $\mu$ is $\sh$-invariant.
\end{theo}
Applying this result also to $a^{-1}$, one deduces that $h_{\mu}(a)=h_{\text{top}}(a)$ if and only if $\mu$ is invariant under the group generated by $\sh$ and $\uh$, and by the unique ergodicity assumption on $G$, $m_{X}$ is the unique measure of that property. It follows that $\mu=m_{X}$. We want to remark at this point the connection to our mixing assumption. The Lie algebra generated by the Lie algebras $\sha$ and $\uha$ of $\sh$ and $\uh$ is a Lie ideal $\gof$ of $\gog=\on{Lie}(G)$, called \textit{Auslander ideal}. If $a$ is mixing on $L^{2}(X)$ then already $\sh$ is uniquely ergodic. On the other hand, the Howe-Moore theorem guarantees mixing in many natural cases (i.e.\ any ergodic action of a simple group is mixing).

We quickly compare the proof of Theorem~\ref{theorem} and Theorem~\ref{pisathm}. The if and only if part of the latter comes from the strict convexity of the map $x\mapsto x\log{x}$. A second order approximation of this function shall enable us to deduce Theorem~\ref{theorem}. More precisely, to compare $\mu$ and $m_{X}$ we relate their conditional measures with respect to a $\sh$-$\it{subordinate}$ $\sigma$-algebra $$\cA=\bigvee_{n=0}^{\infty}a^{-n}\xi$$ on $X$. By $\sh$-subordinate we mean that the atoms of $\cA$ consist of $\sh$-$\it{plaques}$, sets of the form $[x]_{\cA}=xV$ for some $V\subset\sh$ containing a neighborhood of the identity. Here $\xi$ is a countable measurable partition of X that generates the Borel $\sigma$-algebra $\cB$ in the sense that $\bigvee_{n=-\infty}^{\infty}\al[-n]\xi=\cB$ modulo $\mu$.  In our situation $V$ will not depend on $x\in X$. The atoms of $[x]_{\al[-l]\cA}$ will then support more and more of the $\sh$-orbit at $x$ as $l\to\infty$.
By the generating assumption, $$h_{\mu}(a)=H_{\mu}(\cA|\al[-1]\cA)\stackrel{\text{Def}}{=}\int-\log{\mu^{\al[-1]\cA}_{x}([x]_{\cA})}d\mu.$$ 
As for the Haar measure, the inner information function is constant (equal to $h_{m_{X}}(a)$) and this shall allow us to rewrite the entropy difference $h_{m_{X}}(a)-h_{\mu}(a)$ as one integral over $\mu$ (!) of the form $$h_{m_{X}}(a)-h_{\mu}(a)=\int_{X}\sum\log{\frac{q_i}{p_i}}q_i d\mu$$ where $p_{i}$ and $q_{i}$ abbreviate the conditional measures ${(m_X)}^{\al[-1]\cA}_{x}([x_{i}]_{\cA})$ and $\mu^{\al[-1]\cA}_{x}([x_{i}]_{\cA})$, respectively and we sum over $x_{i}$ such that $\bigsqcup [x_{i}]_{\cA}=[x]_{\al[-1]\cA}$.
By convexity, the integrand vanishes if and only if $p_{i}=q_{i}$ for all $i$. An iteration shows that all the conditional measures $\mu_{x}^{{\al[-l]\cA}}$ are Haar, from which one can deduce that $\mu$ is necessarily $\sh$-invariant.

To quantify this we will instead relate the above sum to $\frac{1}{2}\left(\sum|p_{i}-q_{i}|\right)^{2}$ by means of Pinsker's inequality. 
Pushing invariance from the conditional measures to the actual measure will require an exponential rate of equidistribution of the $\sh$-plaques when expanded by conjugation with $a^{-l}$, which explains the restriction to $K$-finite vectors.
We note that in this overview we have made use of a stronger property of $\cA$ that we will actually use: Instead of a $\sigma$-algebra whose atoms are parts of $\sh$-orbits, we will use one whose atoms are contained in the manifold of a larger orbit set of the subgroup $\{g\in G: a^{n}ga^{-n} \text{ stays bounded as } n\to\infty\}$.

Apart from how equidistribution of the test function is achieved, the proof of Theorem~\ref{theorem} will be along the lines of \cite{Polo} and \cite{pisa}, the latter being present in a rather trivialized fashion as we do not have to deal with the boundary of the partition (and therefore handle different null sets with respect to $\mu$, $\mu^{\cA}_{y}$ and $m_{X}$) nor use the fact that ${\al[l]\cA}$ actually generates $\cB$ (for which we pay the price by having $l_{f}$ appear in high degree in the final bound of Theorem 1).

We expect that the result also holds for real quotients. It requires, however, a more careful study of the the corresponding equidistribution theorem. Indeed, in the case discussed here, the profinite structure of $G$ provides us with a decomposition of the space $X$ in orbits of compact groups, which then is no longer given. We note that equidistribution is in fact only needed in $L^2$ and not pointwise (as used here), which will give some leeway to a possible extension.

\section{Notation of and Facts about Linear p-adic Lie Groups}
We quickly introduce the necessary notation and collect some facts about linear p-adic Lie groups along the way. Denote by $\bZ_p$ the p-adic integers of $\bQ_p$ and $|\cdot|_p$ the p-adic norm defined by $d(ax)=|a|_pdx$ if $dx$ is a Haar measure of $\bQ_p$. Equip the vector space $\gom=\Mat_{d}{(\bQ_{p})}$ with the maximum norm $\|X\|:=\max_{i,j}|X_{ij}|_p$ which is again non-Archimedean so that $\|X+Y\|\leq\max{\set{\|X\|,\|Y\|}} \text{ and } \|XY\|\leq\|X\|\|Y\|$. The group of invertible matrices $\GL_{d}{(\bZ_{p})}$ with integral entries is endowed with the bi-invariant metric $d(g,h)=\|g-h\|$ that can be extended to a left-invariant metric on $\GL_{d}{(\bQ_{p})}$ in the following way: Decompose $\GL_{d}{(\bQ_{p})}=\bigsqcup_{n>0} g_{n}\GL_{d}{(\bZ_{p})}$ with $g_{1}=e$ and define 
$$f(g) = \left\{
  \begin{array}{l l}
    d(g,e) & \quad \text{if } g\in \GL_{d}{(\bZ_{p})}\\
    n & \quad \text{if }g\in g_{n}\GL_{d}{(\bZ_{p})} \text{for $n>1$.}
  \end{array}\right. $$
and
$$\tilde{f}(g)=\inf\left\{\sum f(h_{i}): g=h_{1}^{\varepsilon_{1}}\dots h_{k}^{\varepsilon_{k}} \text{ with } \epsilon_{i}\in\{\pm1\},\;h_{i}\in\GL_{d}{(\bQ_{p})}\right\}.$$
Then $d'(g_{1},g_{2})=\tilde{f}(g_{1}^{-1}g_{2})$ defines a proper left-invariant metric that extends $d$. 
We will denote $$K^{\gom}_{k}=p^{k}\Mat_{d}{(\bZ_{p})}$$
and for a $p$-adic Lie group $H=\bH(\bQ_{p})$ where $\bH<\GL_{d}$ is algebraic, the ball of radius $p^{-k}$ in $H$ by $$K^{H}_{k}=(e+K^{\gom}_{k})\cap H$$
and for the Lie algebra $\goh$ of $H$ (defined as the set of tangents of analytic curves in $H$ at the identity) we let $K^{\goh}_{k}=K^{\gom}_{k}\cap\goh$. The exponential function $X\mapsto \exp{X}$ defines a locally analytic isomorphism between the Lie algebra $\goh$ and $H$ for which we refer to \cite{platonov}, p.116 and the references therein. In the non-Archimedean setting $\exp$ turns out to be isometric.
\begin{lemma}
\label{sgpstosgps}
The exponential map $\exp{X}=e+X+\frac{X^2}{2!}+\dots$ is isometric and maps Lie algebras to subgroups. More precisely, for any $k\geq2$ one has that $$\exp{K^{\goh}_{k}}=K^{H}_{k}.$$
\end{lemma}
\proc{Proof.}
The matrix exponential map and logarithm $\log(g)=(g-e)-\frac{(g-e)^2}{2}+\dots$ are mutually inverse and isometric whenever they are defined. Indeed, one calculates that $|n!|_{p}=p^{-\sum_{j\geq1}\lfloor n/p^{j}\rfloor}\geq p^{-{n/(p-1)}}$ and thus by Hadamard's formula the radius of convergence of $\exp$ is at least $p^{-{1/(p-1)}}$ so that $\exp$ converges absolutely on $K^{\gom}_{2}$ (The special case here is $p=2$ where $\exp$ converges only on $K^{\gom}_{2}$, for higher $p$, the domain of convergence is actually $K^{\gom}_{1}$). In particular, for $x\in K^{\gom}_{2}$ we have $\|X^{n}/n!\|\leq p^{n}\|X\|^{n}\leq p^{{-n}}$ and thus $\|\exp(X)-e\|=\max_{n}{\|X^{n}/n!\|}=\|X\|$. Similarly, $\log$ converges on $e+K^{\gom}_{1}$.

As already mentioned, $\exp$ is a local isomorphism between $\goh$ and $H$. First note that this implies $\exp(K^{\goh}_{2})\subset H$. Indeed, assume that $\exp$ and $\log$ are isomorphism of $K^{\goh}_{l}$ and $K^{H}_{l}$ for some sufficiently large $l>0$. Then for any polynomial $f$ that defines $\bH$, and any $X\in K^{\goh}_{2}$, we have $f(\exp{tX})=0$ for all $t\in p^{l}\bZ_{p}$. As $f(\exp{tX})$ is a power series in $t$, this implies $f(\exp{tX})=0$ for all $t\in\bZ_{p}$. On the other hand, if $h\in K^{H}_{2}$ then $\exp{p^{l}\log{h}}=h^{p^{l}}\in K^{H}_{l}$ so that $\log{h}\in p^{-l}\goh=\goh$ and thus $\log h\in K^{H}_{2}$ because the logarithm map is isometric. This shows that $\exp$ is also surjective on $K^{H}_{2}$.

\ep\medbreak

\proc{Remark.}
We want to denote a $\bZ_{p}$-submodule $\goh<K^{\gom}_{2}$ closed under taking commutators a $\bZ_{p}$-Lie algebra. Dynkin's form of the Baker-Campbell-Hausdorff formula for $z=\log{(\exp{x}\exp{y})}$ is $$z=\sum z_{n}, \quad z_{n}=\frac1n\sum_{i+j=n}(x_{i,j}+y_{i,j})$$ with
$$x_{i,j}=\sum \frac{(-1)^{m+1}}{m}\frac{\on{ad}(x)^{i_{1}}\on{ad}(y)^{j_{1}}\dots\on{ad}(x)^{i_{m}}(y)}{i_{1}! j_{1}! \dots i_{m}!}$$
summing over $i_{k}+j_{l}\geq1$ such that $\sum_{k\leq m} i_{k}=i$, $\sum_{l<m} j_{l}=j-1$, $i_{m}\geq1$ and a similar formula holds for $y_{i,j}$ (see p.29 \cite{serre}). Bounding the factorials such as in the previous Lemma implies that $z\in\goh$ so that $\exp{\goh}$ is a $p$-adic Lie group if $\goh$ is a $\bZ_{p}$-Lie algebra.
\medbreak

Denote by $\pi:G\to X$ the projection map to the quotient space and let us scale the Haar measure $m_{G}$ on $G$ to be compatible with the probability measure $m_{X}$, $m_{X}(\pi(B))=m_G(B)$ for any sufficiently small ball $B$ in $G$. We mentioned in the beginning that for a $p$-adic homogeneous space the existence of a finite measure already implies compactness. We give a quick argument (taken from \cite{ratner}) since it also provides an injectivity radius independent of $\Gamma$.
\begin{prop}
\label{compact}
The space $\homq$ is compact with uniform injectivity radius $p^{-2}$.
\end{prop}
\proc{Proof.}
We first claim that for any $\Gamma<G$ the intersection $\Gamma\cap K^{G}_{2}$ must be trivial. As $\Gamma$ is discrete and $K^{G}_{2}$ compact, it suffices to show that $\Gamma\cap K^{G}_{2}$ does not contain any finite non-trivial subgroups. As for any $g\in\Gamma\cap K^{G}_{2}$ we have by denseness of $\bZ$ in $\bZ_{p}$ that $\overline{\{g^{n}\}}=\exp{(\bZ_{p}\log{g})}\subset \Gamma\cap K^{G}_{2}$. Hence any non-trivial subgroup must contain a copy of $\bZ_{p}$.
We may deduce the proven fact also for the conjugated variant, $g^{-1}\Gamma g\cap K^{G}_{2}$, which readily implies that the projection map $\pi$ is injective on the neighbourhood $gK^{G}_{2}$ for any $g\in G$ so that $p^{-2}$ is a uniform injectivity radius. By definition of $m_{X}$, $$m_{X}(\pi(gK^{G}_{2}))=m_G(gK^{G}_{2})\equiv \text{const}.$$ But $m_{X}$ is also finite so that $\hom$ must be covered by only finitely many translates of $K^{G}_{2}$.
\ep\medbreak

\section{Bowen Balls and Topological Entropy $h_{m_{X}}(a)$}
\label{sectionBowenBalls}
Let $\gog=\on{Lie}(G)<\slq$ denote the Lie algebra of $G$. Define the stable horospherical subgroup $$G^{-}_a=\{g\in G \;|\; a^nga^{-n}\to e\text{ as } n\to\infty\}$$ and denote its Lie algebra by $\gog^{-}$. Analogously let $$G^{+}_a=\{g\in G \;|\; a^nga^{-n}\to e\text{ as }n\to-\infty\}$$
define the unstable horospherical subgroup with Lie algebra $\gog^{+}$. Further denote by $G^{0}_{a}$ the group consisting of elements $g\in G$ for which $a^{n}ga^{-n}$ stays bounded for both $n\to\pm\infty$ with Lie algebra $\za$ and finally put
$$P=\set{g\in G \;|\; a^nga^{-n} \text{ stays bounded as } n\to\infty}$$ with Lie algebra $\za+\gog^{-}$ which we will denote by the thickened stable horospherical subgroup. The Lie algebras $\sha$, $\uha$ and $\za$ correspond simply to the eigenspaces of $\on{Ad}_{a}$ of eigenvalues in absolute value smaller than one, bigger than one or equal to one respectively, and thus $\gog=\sha+\za+\uha$. 
For each eigenspace $E_{\lambda}$ of $\on{Ad}_{a}$ there exists a basis of the $\bZ_{p}$-module $E_{\lambda}\cap \Mat_{d}{(\bZ_{p})}$. Let $\{X_{i}\}$ be the union of the resulting basis vectors. We introduce a new family of balls $_{a}K^{G}_{k}$ that are adapted to this eigenbasis,
$$_{a}K^{G}_{k}=\exp{\sum_{i}p^{k}\bZ_{p}X_{i}},$$
where we denote the inner sum also by 
$$_{a}K^{\gog}_{k}=\sum_{i}p^{k}\bZ_{p}X_{i}.$$

\proc{Remarks.}
\label{rmkeigenvalues} 
\begin{enumerate}
\item[(i)]  By choice of the $\bZ_{p}$-eigenbasis the set $_{a}K^{\gog}_{k}$ is a Lie algebra over $\bZ_{p}$ and by the remark below Proposition~\ref{sgpstosgps}, $_{a}K^{G}_{k}$ defines a group. 
\item[(ii)] The smoothness parameter $l_f$ of a function $f$ has been defined regarding the original balls $K^{G}_{k}$. We thus may replace $l_f$ by $l_f+|\nu|$ to ensure that $f$ is $_{a}K_{l_f}$-invariant.
\item[(iii)]  As we will work on with the balls adjusted to $a$ defined above, we will drop the subscript $a$ from $_{a}K$.
\end{enumerate}
\medbreak

Denote $\{u^{-}_{i}\}\in\sha$ those eigenvectors $X_{i}$ that are associated to eigenvalues of absolute values $p^{-\nu_{i}}$ less than one.
We also abuse our previous notation by introducing
\begin{equation}\label{notation}K^{\sha}_{k+\nu}=\sum p^{k+\nu_{i}}\bZ_{p}u^{-}_{i}\end{equation}
for the vector $\nu=(\nu_{1},\dots,\nu_{\on{dim}(\sha)})$ with norm $|\nu|=\sum\nu_{i}$. These boxes coincide with the image of $K^{\sha}_k$ under $\on{Ad}_a$. Similarly define $K^{\uha}_{k}$ and $K^{P}_{k,l}=\exp{(K^{\za}_{k}+{K}^{\sha}_{l})}$ with image under $\on{Ad}_a$ equal to $K^{P}_{k,l+\nu}$ using the eigenvectors from above.
We introduced the new basis of the Lie algebra so that we are able to split the ball ${K}^G_{k}$ into a thickened stable and unstable component as such:
\begin{lemma}
\label{balldecomposition} The decomposition
$K^{\uh}_{k} {K}^P_{k,k}= {K}^G_{k}$ holds for all $k>1$.
\end{lemma}
\proc{Proof.}
We only have to address the inclusion ${K}^G_{k}\subset {K}^{\uh}_{k} {K}^P_{k,k}$. For this let $$g=\exp{(v_{0}+w_{0})}\in{K}^G_{k}$$ where $v_{0}\in\sha$ and $w_{0}\in\uha+\za$, both of norm $\leq p^{-k}$. Define $f_{0}=\exp{v_{0}}$ and $h_{0}=\exp{w_{0}}$ then $f_{0}^{-1}gh_{0}^{-1}=\exp{(v_{1}+w_{1})}$ for some $v_{1}\in\sha$ and $w_{1}\in\uha+\za$ with norm at most $p^{-2k}$ by application of the Baker-Campbell-Hausdorff formula mentioned in the remark below Proposition~\ref{sgpstosgps}.
Continuing this procedure we find sequences $f_{i}\in {K}^{\uh}_{k} $ and $h_{i}\in {K}^P_{k,k}$ such that $$f_{i}^{-1}\dots f_{0}^{-1}gh_{0}^{-1}\dots h_{i}^{-1}\to e.$$ On the other hand, $F_{i}=f_{0}\dots f_{i}$ and $H_{i}=h_{0}\dots h_{i}$ all lie in the compact set ${K}^G_{k}$ so that we find a converging subsequence $F_{i}\to f\in {K}^{\uh}_{k}$ and $H_{i}\to h\in {K}^P_{k,k}$ for which also $f^{-1}gh^{-1}=e$ holds.
\ep\medbreak

Let $\cP^X_2$ denote the partition of $X$ into $m_{G}(K^{G}_{2})^{-1}$ balls and $\cP^X_k$ for $k>2$ the corresponding refinement in which each atom of $\cP^X_2$ is split into $p^{(k-2)d}$ smaller balls where $d$ is the dimension of $G$. This is indeed well formulated: $K^{\gog}_{k}$ is the disjoint union of $p^{d}$ translates of $K^{\gog}_{k+1}$ (with distance $p^{-k}$ from each other) so that after applying the exponential map, $K^{G}_{k}$ splits into $p^{d}$ copies of $K^{G}_{k+1}$ and after projecting, being inside the injectivity radius, the property that any two balls either coincide or are disjoint is passed onto atoms of $\cP^X_k$ and thus gives a unique refinement. 

Define the $\sigma$-algebra generated by $\cP_{k}^{X}$ under $a$ $$\cA=(\cP^X_k)^\infty_0=\bigvee^\infty_{l=0}\al\pk$$ to be the smallest $\sigma$-algebra containing all the partitions $\al\pk$. If $[x]_{\pk}=x{K}^G_{k}$ denotes the atom of $\pk$ then 
\begin{equation}
\label{eq:atom}
[x]_{\al\pk}=\al{}[\al[l]x]_{\pk}=xa^{-l}{K}^G_{k}a^{l}.
\end{equation}
In order to understand the form of an atom $[x]_{\cA}$ we need the following lemma.
\begin{lemma} For $k\in\bZ_{>1}$ such that $k-\max{\nu_{i}}>1$, it holds for all $n\in\bZ_{\geq0}$ that
$$[x]_{(\pk)_{0}^{n}}=x\exp{\left({K}^{\uha}_{k+n\nu}+{K}^{\za}_{k}+{K}^{\sha}_{k}\right)}.$$
\end{lemma}
\proc{Proof.}
We start by noting that $$\exp{\left({K}^{\uha}_{k+n\nu}+{K}^{\za}_{k}+{K}^{\sha}_{k}\right)}$$ is the intersection of $a^{-l}{K}^G_{k}a^{l}$ for all $0\leq l\leq n$. Indeed, by our notation introduced in line (\ref{notation}) we have
\begin{align*}\bigcap_{l\leq n}a^{-l}\left({K}^{\uha}_{k}+{K}^{\za}_{k}+{K}^{\sha}_{k}\right)a^{l}&= \bigcap_{l\leq n}\left({K}^{\uha}_{k+l\nu}+{K}^{\za}_{k}+{K}^{\sha}_{k-l\nu}\right)\\&={K}^{\uha}_{k+n\nu}+{K}^{\za}_{k}+{K}^{\sha}_{k}
\end{align*}
so that
$$\bigcap_{l\leq n}a^{-l}\exp{\left({K}^{\uha}_{k}+{K}^{\za}_{k}+{K}^{\sha}_{k}\right)}a^{l}=\exp{\left({K}^{\uha}_{k+n\nu}+{K}^{\za}_{k}+{K}^{\sha}_{k}\right)}.$$
We only have to take care of the inclusion $\subset$ in the statement of the lemma. For this  let $y\in[x]_{(\pk)_{0}^{n}}$ which means by definition that for all $l\leq n$ the points $\al[l]x,\al[l]y$ lie in the same partition element of $\pk$. In particular, if $y=xg_{0}$ for $\|g_{0}\|\leq p^{-k}$ then also $\al[{}]y$ and $\al[{}]x$ are $p^{-k}$-close, i.e. there exists $g_{1}\in{K}^{G}_{k}$ such that $\al[{}]y=xg_{0}a^{-1}=xa^{-1}g_{1}$, or equivalently,
$$xg_{0}=xa^{-1}g_{1}a.$$
We chose $k$ such that $a^{-1}g_{1}a\in K^{G}_{2}$. Using the bound on the injectivity radius from Proposition \ref{compact} the above equation lifts to the group level
$g_{0}=a^{-1}g_{1}a$
from which we conclude that $$g_{0}\in {K}^{G}_{k}\cap a^{-1}{K}^{G}_{k}a.$$
Repeating this step for the points $y_{1}=x_{1}g_{1}$ and $x_{1}=xa^{-1}$ we find $g_{2}\in{K}^{G}_{k}$ such that $y_{1}a^{-1}=x_{1}a^{-1}g_{2}$, or again equivalently that $g_{1}\in a^{-1}{K}^{G}_{k}a$ and thus $$g_{0}\in {K}^{G}_{k}\cap a^{-1}{K}^{G}_{k}a\cap a^{-2}{K}^{G}_{k}a^{2}.$$
Continuing this argument we find successively $g_{0},\dots, g_{n}$ all elements of ${K}^{G}_{k}$ with relation $g_{i}=a^{-1}g_{i+1}a$ so that $g_{0}$ lies indeed in $\exp{\left({K}^{\uha}_{k+n\nu}+{K}^{\za}_{k}+{K}^{\sha}_{k}\right)}$.
\ep\medbreak

Letting $n\to\infty$ we see that the atoms $$[x]_{\cA}=x\exp{\left({K}^{\za}_{k}+{K}^{\sha}_{k}\right)}=x{K}^{P}_{k}$$ are balls in the thickened horospherical direction and for $k$ large enough (depending on $\nu$),
\begin{eqnarray}
\label{eq:adjointaction}
[x]_{\al[-1]\cA}&=x\exp{\left(K_{k}^{\za}+\sum p^{k-\nu_{i}}\bZ_{p}u_{i}\right)}=xK^P_{k,k-\nu}=\bigsqcup_{j=1}^{p^{|\nu|}}xg_jK^P_{k,k}
\end{eqnarray}
for some $p^{|\nu|}$ elements $g_j\in\sh$.

We use Bowen's formalism for \text{homogeneous} measures \cite{bowen} to calculate that $$h_{m_{X}}(a)=|\nu|\log{p}$$ where $p^{|\nu|}$ is the product of the absolute values of all eigenvalues of $\on{Ad}_{a}$ (with multiplicities) that are greater than one.
\begin{prop}
\label{bowenentropy}We have $h_{m_{X}}(a)=h_{m_{\gog}}(\on{Ad}_{a})=|\nu|\log{p}$.
\end{prop}
\proc{Proof.}
As the space $X$ is compact and since the projection map $\pi$ is locally isometric by the construction of the metric on $X$, $d_{X}(\Gamma g,\Gamma h)=\inf_{\gamma\in\Gamma}d(g,\gamma h)$ we might as well calculate the entropy of $a$ acting on $G$. If $D_{n}(e,k,a)=\bigcap_{l=0}^{n-1}a^{-l}K^G_{k}a^{l}$ denotes a family of Bowen balls in $G$, the topological entropy is $$\lim_{k\to\infty}\limsup_{n\to\infty}-\frac{1}{n}\log{m_{G}(D_{n}(e,k,a))}.$$
From above we have $D_{n}(e,k,a)=\exp{\left(K_{k+(n-1)\nu}^{\uha}+K_{k}^{\za}+K_{k}^{\sha}\right)}$.
Thus, $D_{0}(e,k,a)$ is the disjoint union of $p^{(n-1)|\nu|}$ translates of the $n$th Bowen ball so that $m_{G}(D_{n}(e,k,a))=p^{-(n-1)|\nu|}m_{G}(D_{0}(e,k,a)),$
concluding that 
\begin{eqnarray*}
\lim_{n\to\infty}-\frac{1}{n}\log{m_{G}(D_{n}(e,k,a))}&=&\lim_{n\to\infty}\left(\frac{n-1}{n}|\nu|\log{p}-\frac{1}{n}\log{m_{G}(D_{0}(e,k,a))}\right)\\&=&|\nu|\log{p}-0.
\end{eqnarray*}
\ep\medbreak
This calculation shows that the modular function $\on{mod}(a,\sh)$ of $a$ with respect to the inner action on $\sh$ is $p^{|\nu|}$ so that $$h_{\text{top}}(a)=h_{m_{X}}(a)=\log{\operatorname{mod}(a,\sh)}$$ as claimed earlier.

\section{Entropy generating Partition}
\label{generator}
To calculate the entropy one usually finds a suitable generating partition. From the preceding section we see that atoms of the $\sigma$-algebra $\bigvee_{n=-\infty}^\infty \al[-n]\pk$ are \textit{plaques} of the form $xK_k^{G^{0}_a}$. In particular, $\pk$ does not generate the Borel $\sigma$-algebra under the action of $a$. However, the following holds
\begin{prop}
\label{generatorprop}
A fine enough partition $\pk$ will still be entropy-generating for any $a$-invariant measure $\mu$, in the sense that
$$h_\mu(a)=h_\mu(a,\pk)=H_\mu\left(\pk|\textstyle\bigvee_{n=1}^\infty \al[-n]\pk\right).$$
\end{prop}
Any $k\geq|\nu|+2$ works and thanks to the concrete description of the plaques, this is shown rather painlessly (compare to the real analogue, Proposition 9.2 in \cite{EKL}). 
\proc{Proof.}
First take the increasing sequence of $\sigma$-algebras $\sigma(\cP^X_l)$ that converges to the Borel $\sigma$-algebra of $X$ as $l\to\infty$ so that $h_\mu(a)=\lim_{l\to\infty}h_\mu(a,\cP^X_l).$
On the other hand,
$$h_\mu(a,\cP^X_l)\leq h_\mu(a,\cP^X_l\vee \cP^X_k)=h_\mu(a,\cP^X_k)+h_\mu\left(a,\cP^X_l|\textstyle\bigvee_{n=-\infty}^\infty \al[-n]\pk\right)$$
and we claim the latter term vanishes which then by definition of the $h$-entropy implies that $h_\mu(a)=h_\mu(a,\pk)$. Indeed,
$$h_\mu\left(a,\cP^X_l|\textstyle\bigvee_{n=-\infty}^\infty \al[-n]\pk\right)=H_\mu\left(\cP^X_l|\textstyle\bigvee_{m=1}^\infty\al[-m]\cP^X_l\vee\textstyle\bigvee_{n=-\infty}^\infty \al[-n]\pk\right)$$
is zero if the partition in the first argument of $H_\mu$ is contained in the $\sigma$-algebra of the second argument. But this is an easy calculation for our explicit description of the atoms,
$$[x]_{\bigvee_{m=1}^\infty\al[-m]\cP^X_l}\cap[x]_{\bigvee_{n=-\infty}^\infty \al[-n]\pk}=xK^P_{l,l-\nu}\cap xK^{G^{0}_a}_k=xK^{G^{0}_a}_l\subset[x]_{\cP^X_l}$$
for any $l\geq k$ and $l\geq|\nu|+2$ by equation $(\ref{eq:adjointaction})$.
\ep\medbreak
\section{Proof}
We now begin the proof of Theorem~\ref{theorem}. We will first rewrite the integral of $f$ over the Haar measure of $X$ as limit of the conditional expectations (whose properties we recall below) with respect to the family of $\sigma$-algebras $\al[-n]\cA$,
$$f_{n}(x):=E_{m_{X}}(f|\al[-n]\cA)(x).$$ Note also that we will ultimately assume that $m_{X}(f)=0$.
\label{proof}
\subsection{Contracting to (almost) $\sh$-orbit integrals.}
\label{orbitintegrals}
We will now bound the expression $\left|m_{X}(f)-\mu(f)\right|$ by looking at how $f$ integrated over thickened $\sh$-orbits behaves with respect to $\mu$ on average. For that, fix $k=2+|\nu|+l_{f}$ forcing $f$ to be $\cA$-measurable and the exponential map to be defined on atoms of $\cA=\left(\pk\right)^\infty_0$. We recall that for a measure $\nu$ and a $\sigma$-algebra $\cC$ on a measure space $Y$, the family of conditional measures $\set{\nu^{\cC}_{y}}_{y\in Y}$ is uniquely determined by the properties of the conditional expectation $E_{\nu}(f|\cC)$ defined by
$$y\mapsto \int f\;d\nu^{\cC}_{y}$$
to exist $\nu$-a.e. and to be $\cC$-measurable and integrable for every $\nu$-integrable $f$, and that the integral equation
$$\int_{C}\int f\;d\nu^{\cC}_{y}d\nu(y)=\int_{C}f\;d\nu$$
holds for every $C\in\cC$.
\begin{lemma} For an $\cA$-integrable function $f$ and an $a$-invariant measure $\mu$ it holds that
$$\left|m_{X}(f)-\mu(f)\right|\leq\sum_{n\geq0}\int_{X} \left|f_{n+1}(x)-f_{n}(x)\right|d\mu$$
with $f_{n}$ as defined above.
\end{lemma}
\proc{Proof.}
We have seen that the atoms $[x]_{\cA}$ of $\cA$ are $xK^{P}_{k}$. Because $K^{P}_{k}$ is a subgroup this implies that the conditional measure $(m_{X})^{\cA}_{x}$ is just the push forward of the map $g\mapsto g\acts x$ from $P$ to $X$ of the Haar measure $m_{P}$ restricted to $K^{P}_{k}$,
$$(m_{X})^{\cA}_{x}=\frac{1}{m_{P}(\Pk)}m_{P}\acts x.$$
Indeed, the function $\frac{1}{m_{P}(\Pk)}\int_{\Pk}f(g\acts x)dm_{P}(g)$ is constant under the action of $K^{P}_{k}$ and thus $\cA$-measurable and by an application of Fubini and invariance of $m_{G}$ for every measurable $A\acts y=B\Pk\acts y\in\cA$, $B\subset K^{\uh}_k$ (which generate $\cA$),
\begin{multline*}
\frac{1}{m_{P}(\Pk)}\int_{A\acts y}\int_{\Pk}f(g\acts x)dm_{P}(g)dm_{X}(x)\\=\frac{1}{m_{P}(\Pk)}\int_{\Pk}\int_{A}f(yh^{-1}g^{-1})dm_{G}(h)dm_{P}(g) \\=\int_{A}f(h\acts y)dm_{G}(h)
\end{multline*}
so that $m_{P}\acts x$ satisfies the two properties that uniquely characterize the conditional measure. 
Note that $m_P$ is unimodular when restricted to the compact plaques $a^{-n}K^P_ka^n$. Now $(m_{X})_{x}^{\al[-l]\cA}=\al[-l](m_{X})^{\cA}_{\al[l]x}$ and we similarly check that
\begin{eqnarray*}
f_{n}(x)=E_{m_{X}}(f|\al[-n]\cA)(x)&=&\frac{1}{m_{P}(\Pk)}\int_{\Pk}f(\al[-n](g\acts(\al[n]x))dm_{P}\\&=&\frac{1}{m_{P}(a^{-n}\Pk a^{n})}\int_{a^{-n}\Pk a^{n}}f(g\acts x)dm_{P}.
\end{eqnarray*}

We will understand the last expression as a thickened horospherical flow and show in Theorem~\ref{effequ} that it equidistributes uniformly, i.e.\ $f_{n}(x)\to m_{X}(f)$ as $n\to\infty$ independent of $x$. Since $f$ is $\cA$-measurable, $f_0=f$. Write 
$$\left|m_{X}(f)-\mu(f)\right|=\left|\mu(m_{X}(f)-f)\right|=\lim_{k\to\infty}\left|\mu(f_{k}-f)\right|$$
and since
$$f_k-f=\sum_{n=0}^{k-1}(f_{n+1}-f_n)+f_0-f=\sum_{n=0}^{k-1}(f_{n+1}-f_n),$$
we see that we indeed moved the problem to understanding the sum of the  differences $\int_{X}\left|f_{n+1}(x)-f_{n}(x)\right|d\mu$.
\ep\medbreak

\subsection{Bounding the term $f_{n+1}-f_{n}$.}
\label{telescoping}
We make use of the following relation:
\begin{lemma}
If $f\circ a^{n}$ denotes the map $x\mapsto f(\al[-n]x)$ then the recursion formula
$$f_{n+1}(x)=E_{m_X}(f_{n}\circ a^{n}|\al[-1]\cA)(\al[n]x)$$
holds for all $n\geq0$.
\end{lemma}
\proc{Proof.}
The right-hand side is equal to
\begin{eqnarray*}
\int_{a^{-1}\Pk a}f_{n}(a^{-n}\acts (g\acts (\al[n]x)))dm_{P}&=&\int_{a^{{-(n+1)}}\Pk a^{{n+1}}}f_{n}(g\acts x)dm_{P}\\
&=&E_{m_X}(f_{n}|a^{(n+1)}\acts\cA)(x)=E_{m_X}(f|a^{(n+1)}\acts\cA)(x).
\end{eqnarray*}
\ep\medbreak
The following classical lemma in information theory plays the crucial role to get a quantitative estimate from the convexity of the entropy function.
\begin{lemma}[Pinsker's Inequality, Lemma 12.6.1 \cite{cover}]
\label{lem:convexity}
Let $$\phi_{\underline{p}}(\underline{q})=-\sum\log{\frac{p_i}{q_i}}q_i$$ be defined on the $n-1$-dimensional simplex of probability vectors $\underline{q}=(q_1,\dots,q_{n})$ such that the $q_i$'s are positive and sum up to one. Then for all probability vectors $\underline{q}$ and $\underline{p}$ it holds $$\|\underline{p}-\underline{q}\|_{1}^2\leq 2\phi_{\underline{p}}(\underline{q}).$$
\end{lemma}
We will apply this to the probability vectors defined by the coordinates $$p_i(x)={(m_X)}^{\al[-1]\cA}_{x}([x_{i}]_{\cA}) \quad\text{and} \quad\quad q_i(x)=\mu^{\al[-1]\cA}_{x}([x_{i}]_{\cA}).$$
\begin{lemma}
\label{lem:entropybound}
With assumptions on $f$ and $\mu$ as before and $\phi_{p}$ defined as in Lemma $\ref{lem:convexity}$ it holds for all $n\geq0$ that
$$\int_{X} \left|f_{n+1}(x)-f_{n}(x)\right|d\mu\leq\sqrt{2}\|f_{n}\|_\infty\left(\int_{X}\phi_{\underline{p}(x)}(\underline{q}(x))dm_{X}\right)^{\frac12}.$$
\end{lemma}
\proc{Proof.}
With the new expression of $f_{n+1}$, we write
\begin{align*}
\int_{X} \left(f_{n+1}(x)-f_{n}(x)\right)d\mu=&\int_{X} \left(f_{n+1}(\al[-n]x)-f_{n}(\al[-n]x)\right)d\mu\\=&\int_{X} (E_{m_{X}}(f_{n}\circ a^{n}|\al[-1]\cA)(x)-E_{\mu}(f_{n}\circ a^{n}|\al[-1]\cA)(x))\;d\mu.
\end{align*}
In order to compare the conditional expectations of the $\cA$-measurable function $f_{n}\circ a^{n}$ with respect to $\mu$ and $m_{X}$, we decompose
$$E_{m_X}(f_{n}\circ a^{n}|\al[-1]\cA)(x)=\sum_{j} f_{n}(\al[-n]x_{j}){(m_X)}^{\al[-1]\cA}_{x}([x_{j}]_{\cA})$$
and
$$E_{\mu}(f_{n}\circ a^{n}|\al[-1]\cA)(x)=\sum_{j} f_{n}(\al[-n]x_{j})\mu^{\al[-1]\cA}_{x}([x_{j}]_{\cA})$$
where $x_j=xg_j\in \hom$ are the $p^{|v|}$ points from $(\ref{eq:adjointaction})$ to represent the atoms $[x_{j}]_{\cA}$ of the partitioning $[x]_{\al[-1]\cA}$ in $\cA$.
By Lemma \ref{lem:convexity},
\begin{multline}
\label{bootstrap}
\int_{X} \left|f_{n+1}(x)-f_{n}(x)\right|d\mu\\
\leq\int_{X}\sum_{j} \left|f_{n}(\al[-n]x_{j})\right|\left|{(m_X)}^{\al[-1]\cA}_{x}([x_{j}]_{\cA})-\mu^{\al[-1]\cA}_{x}([x_{j}]_{\cA})\right|dm_{X}\\
\leq \sqrt{2}\|f_{n}\|_\infty\int_{X}\sqrt{\phi_{\underline{p}(x)}(\underline{q}(x))}dm_{X}.
\end{multline}
Finally, $\int_{X}(\phi_{\underline{p}(x)}(\underline{q}(x)))^{\frac12}dm_{X}\leq\left(\int_{X}\phi_{\underline{p}(x)}(\underline{q}(x))dm_{X}\right)^{\frac12}$ by the Cauchy-Schwarz inequality and $m_{X}(X)=1$.
\ep\medbreak
\subsection{Relating $f_{n+1}-f_{n}$ to the entropy difference.}
The relative entropy $\phi_{\underline{p}}(\underline{q})$ of the two distributions $\underline{p}$ and $\underline{q}$ relates to their entropies as follows.
\begin{lemma} For the $\sigma$-algebra $\cA=\left(\pk\right)^\infty_0$ constructed from the entropy generating partition $\pk$  (see Section \ref{generator}) the following equality holds
$$h_{m_{X}}(a)-h_\mu(a)=\int\phi_{\underline{p}(x)}(\underline{q}(x))d\mu(x).$$
\end{lemma}
\proc{Proof.}
By assumption on the entropy generation we can write the entropy $h_\mu(a)$ as
\begin{align*}
H_\mu(\cA|\al[-1]\cA)=&\int-\log \mu^{\al[-1]\cA}_{x}([x]_{\cA})d\mu(x)\\=&\int\int-\log \mu^{\al[-1]\cA}_{y}([y]_{\cA})d\mu^{\al[-1]\cA}_{x}(y)\;d\mu(x)\\
=&\int\int-\log \mu^{\al[-1]\cA}_{x}([y]_{\cA})d\mu^{\al[-1]\cA}_{x}(y)\;d\mu(x)\\
=&\int\sum_{i=1}^{p^{|\nu|}}\left(-\log \mu^{\al[-1]\cA}_{x}([x_{i}]_{\cA})\mu^{\al[-1]\cA}_{x}([x_i]_{\cA})\right)d\mu(x),
\end{align*}
where the first equality sign is just the definition of $H$, the third equality follows from the fact that the conditional measures $\mu^{\al[-1]\cA}_{x}$ has support on $[x]_{\al[-1]\cA}$ and for almost every point there is only one conditional measure on a particular atom, i.e.\ $\mu^{\al[-1]\cA}_{x}=\mu^{\al[-1]\cA}_{y}$ if $[x]_{\al[-1]\cA}=[y]_{\al[-1]\cA}$ for $\mu-$a.e. $x,y\in X$. The final equality follows now from the fact that $-\log \mu^{\al[-1]\cA}_{x}([y]_{\cA})$ is constant on $\cA$ atoms and $[x]_{\al[-1]\cA}$ decomposes exactly into $p^{|\nu|}$ of those.

We repeat the last calculation for the Haar measure $m_{X}$ and note that by Section \ref{orbitintegrals} the information function with respect to the Haar measure $\log {(m_{X})}^{\al[-1]\cA}_{x}([y]_{\cA})$ is constant (and equal to $|\nu|\log{p}$) so that we may integrate this information function also against $\mu$ instead of $m_{X}$. We therefore may write the difference of the entropies in terms of $\phi$:
\begin{eqnarray*}
h_{m_{X}}(a)-h_\mu(a)&=&\int -\log {(m_X)}^{\al[-1]\cA}_{x}([x]_{\cA})+\log\mu^{\al[-1]\cA}_{x}([x]_{\cA}) d\mu(x)\\
&=&\int\sum_i-\log\frac{{(m_X)}^{\al[-1]\cA}_{x}([x_i]_{\cA})}{\mu^{\al[-1]\cA}_{x}([x_i]_{\cA})}\mu^{\al[-1]\cA}_{x}([x_i]_{\cA})d\mu(x)\\
&=&\int\phi_{\underline{p}(x)}(\underline{q}(x))d\mu(x).
\end{eqnarray*}
\ep\medbreak

We obtain therefore the bound
\begin{equation}
\label{concludingbound}
\left|m_{X}(f)-\mu(f)\right|\leq\sum_{n=0}^{\infty}\left|\mu(f_{n+1}-f_n)\right|\leq\sqrt{2} \sum_{n=0}^{\infty}\|f_{n}\|_\infty\left(h_{m_{X}}(a)-h_{\mu}(a)\right)^{\frac{1}{2}}
\end{equation}
which is only useful if one can show that $$f_n(x)=\frac{1}{m_P(a^{-n}K^P_{k,k}a^n)}\int_{a^{-n}K^P_{k,k}a^n}f(g\acts x)dm_{P}(g)$$ decays fast enough to zero (namely so that $\sum_{n=0}^{\infty}\|f_{n}\|_\infty<\infty)$. Since we may assume that $m_{X}(f)=0$, as Theorem~\ref{theorem} is trivial for constant functions, it remains to prove effective equidistribution along the thickened horospherical orbits for the functions $f\in L^2_0(X)$ of vanishing integral. As we will see, the rate of convergence will depend on $\Gamma$.

\subsection{Effective equidistribution.}
It remains therefore to prove the following equidistribution theorem.
\begin{theo}
\label{effequ} Assume that the action of $G$ on $L^{2}(X)$ is exponentially fast mixing, that is, there exist strictly positive constants $c$,$\alpha$ and $\delta$ such that for any locally constant $f,h\in L^{2}_{0}(X)$ with degree of smoothness $l_{f}$ resp.\ $l_{h}$ and for all $n\in\bZ$ one has $$|\langle f\circ a^n,h\rangle|\leq c p^{(l_{f}+l_{h})\alpha}\|f\|_{L^{2}(X)}\|h\|_{L^{2}(X)}\|a\|^{-\delta n}.$$
Then, with notation as before, for $k= l_f+|\nu|+2$ we have effective equidistribution of the family of sets $a^{-n}K^P_{k,k}a^n$,
\begin{multline*}\left\|\frac{1}{m_P(a^{-n}K^P_{k,k}a^n)}\int_{a^{-n}K^P_{k,k}a^n}f(g\acts x)dm_{P}(g)\right\|_\infty\\\leq\frac{c}{\sqrt{m_G(K^G_{2})}}p^{(\alpha+\frac d2)|\nu|+2\alpha}p^{l_{f}(2\alpha+\frac{d}{2})}\|f\|_{L^{2}(X)}\|a\|^{-\delta n}.\end{multline*}
\end{theo}
We apply Margulis' trick used to prove equidistribution of the horospherical flow via mixing. A presentation of this method can be found in \cite{vol1}, Chapter 11.
\begin{prop}
\label{margulistrick}
The left-hand side of the claimed inequality in Theorem~\ref{effequ} is equal to the sup-norm of the matrix coefficient $\langle f\circ a^{n},h_{x}\rangle_{L^{2}(X,m_{X})}$, where $h_{x}$ is a smooth function depending on $\al[n] x$.
\end{prop}
The dependencies on $h_{x}$ shall not disconcert us further as the force of the mixing assumption lies in the fact that the bounds given only depend on the norm of $f$ and $h_{x}$.
\proc{Proof.}[Proof of Proposition~\ref{margulistrick}]
We may distort the function by $g^{+}$ for any $g^{+}\in a^{{-n}}K^{\uh}_{k}a^n$, as $k$ is assumed to be larger then the degree of smoothness of $f$ and conjugation by $a^{{-n}}$ only abates the amount of distortion by $\uh$ for positive $n$, the conditional expectation $f_n$ at $x$ equals
$$\frac{1}{m_{\uh}(a^{{-n}}K^{\uh}_{k}a^n)m_P(a^{-n}K^P_{k,k}a^n)}\int_{a^{-n}K^{P}_{k,k}a^n}\int_{a^{{-n}}K^{\uh}_{k}a^n}f(g^{+}g\acts x)dm_{\uh}(g^{+})dm_{P}(g).$$
We now apply the following lemma (found e.g.\ as Lemma 11.31 in \cite{vol1}) to rewrite the product $dm_{\uh}dm_{P}$ as $dm_{G}$.
\begin{lemma}
For any closed subgroups $S$ and $T$ in $G$ such that $S\cap T=\set{e}$ and $ST$ contains a neighborhood of $e$, the Haar measure $m_{G}$ restricted to $ST$ is proportional to the push forward of $m^{l}_{S}\times m^{r}_{T}$ under the product map $(s,t)\mapsto st$ where $m^{l}_{S}$ is the left Haar measure on $S$ and $m^{r}_{T}$ is the right Haar measure.
\end{lemma}
We set $S=\uh$ which is unimodular and $T=P$ for which we had chosen the right Haar measure in Section \ref{orbitintegrals} and thus
$$f_n(x)=\frac{1}{m_{G}(a^{{-n}}K^{\uh}_{k}K^P_{k,k}a^n)}\int_{a^{{-n}}K^{\uh}_{k}K^P_{k,k}a^n}f(g\acts x)dm_{G}(g).$$
By Lemma \ref{balldecomposition} we have that $K^{\uh}_{k}K^P_{k,k}=K^G_{k}$. We find that $f_{n}$ is equal to
$$\frac{1}{m_{G}(K^G_{k})}\int_{X}f(a^{{-n}}\acts y)\mathbbm{1}_{K^G_{k}a^{n}\acts x}(y)dm_{X}(y).$$
This is a matrix coefficient under the regular representation at $a^{-n}$ with smooth vectors $f$ and $h_{x}$ defined by $h_{x}(y)=\frac{1}{m_{G}(K^G_{k})}\mathbbm{1}_{K^G_{k}a^{n}\acts x}(y)$.
\ep\medbreak
\proc{Proof.}[Proof of Theorem~\ref{effequ}]
Recall that we assume that $m_G$ is compatible with $m_X$, i.e. $m_G(K^G_{k})=m_X(K^G_{k}\acts x)$. Further $h_{x}-1\in L^{2}_{0}$ and is ramified of level $k$ , $\langle f\circ a^{n},h_{x}\rangle_{L^{2}(X,m_{X})}=\langle f\circ a^{n},h_{x}-1\rangle_{L^{2}(X,m_{X})}$ and since $m_G(K^G_{k})=p^{-d(l_f+|\nu|)}m_G(K^G_{2})$ we find that $$\|h_{x}-1\|_{L^{2}(X)}=\sqrt{\frac{1}{m_G(K^G_{k})}-1}\leq\frac{1}{\sqrt{m_G(K^G_{2})}}p^{d(l_f+|\nu|)/2}$$
and consequentially,
$$\left|f_{n}(x)\right|\leq\frac{c}{\sqrt{m_G(K^G_{2})}}p^{(\alpha+\frac d2)|\nu|+2\alpha}p^{l_{f}(2\alpha+\frac{d}{2})}\|f\|_{L^{2}(X)}\|a\|^{-\delta n}$$
which proves Theorem~\ref{effequ}.
\ep\medbreak
\proc{Proof.}[Proof of Theorem~\ref{theorem}]
In particular, we can conclude Theorem~\ref{theorem} since we found a bound for the sum appearing in line (\ref{concludingbound}), namely $$\sum\left\|f_{n}\right\|\leq\frac{c}{\sqrt{m_G(K^G_{2})}}\frac{1}{1-\|a\|^{\delta}}p^{(\alpha+\frac d2)|\nu|+2\alpha}p^{l_{f}(2\alpha+\frac{d}{2})}\|f\|_{L^{2}(X)}.$$ As mentioned in the remark of Section~\ref{sectionBowenBalls} we adjusted the balls $K^G_k$ to $a$. To reflect this change we replaced (implicitely) $l_f$ by $l_f+|\nu|$. This and the fact that $|\nu|\log p=h_{m_{X}}$ gives the expression for $\kappa$ as in the statement of the theorem.
\ep\medbreak

\section{Effective Decay of Matrix Coefficients for $\on{SL}_{m}(\bQ_{p})$}

We finish by giving some references for the exponentially fast mixing property of unitary representations of simple groups and particularly consider $\on{SL}_{m}(\bQ_{p})$. In Howe and Moore (\cite{hm}) it is shown that the matrix coefficients of a unitary representation of a simple algebraic group (over any local field) without nonzero invariant vectors vanish. For non-Archimedean fields Lemma X.3.4, \cite{wallach} is used: For any irreducible smooth and admissible representation $\rho$ there exists $t>0$ such that for any vectors $v$ and $w$
$$|\langle\rho(g)v,w\rangle|\ll\|v\|\|w\|\Xi(g)^t$$
for any $g$. The assumptions that $\rho$ is smooth and admissible mean that every vector $v$ is smooth, and the dimension of $K_{l}v$ is finite for any $l$. The bounding function $\Xi$ is the Harish-Chandra function, defined by
$$\Xi(g)=\int_K\delta_P^{1/2}(gk)dk$$
where in case of $G=\SLQ$, the function $\delta_P$ is the modular character of the group of upper triangle matrices $P$. By Theorem 4.2.1 in \cite{silberger} it holds that for any  $\alpha>-\frac12$  $$\Xi(g)\ll_\alpha\|g\|^{\alpha}.$$
In fact, the regular representation of any locally compact group $G$ with co-compact lattice has a spectral gap in this sense (\cite{margulisbook}, Chapter III.1). Note that in the form of the statement, the implicit constant will depend on the rate of smoothness of $v$ and $w$. 
This can however be made explicit, for example if the representation is the Koopman representation as in our setup, see Proposition 2.5 in \cite{shalom} or more generally, as in the argument for Theorem 2 in \cite{chh}.
In order to get hold of $\alpha$ we want to cite the following result for higher-rank groups $\SLQ$, $m>2$, which enjoy Property (T) so that one can give not only an explicit but also a uniform bound (with respect to various $\Gamma$). 
\begin{theo}
\label{decay3} Let $m>2$. For any smooth functions $f$ and $g$ in $L^2_0(X)$ both fixed by $K^G_l$ it holds for any non-negative integer $n$ that
$$|\langle f\circ a^n,g\rangle\ll p^{l(m^2-1)}\|f\|_{L^{2}(X)}\|g\|_{L^{2}(X)}\|a\|^{-\delta n}$$
for some (explicit) $\delta>0$, where the implicit constant and $\delta$ are independent of $\Gamma$.
\end{theo}
\proc{Proof.}
A bound for matrix coefficients valid uniformly among unitary representations for higher rank (real) semi-simple groups is well known and is grounded on property (T) of Kazhdan and observed for example by Cowling, \cite{cowling}. This has been generalized to reductive groups over arbitrary local fields by Oh \cite{oh}: For any unitary representation $\rho$ of $\SLQ$ without invariant vector, and any $K$-finite unit vectors $v$ and $w$ of $\rho$,
$$|\langle\rho(g)v,w\rangle|\leq(\on{dim}{Kv})^\frac{1}{2}(\on{dim}{Kw})^\frac{1}{2}\prod_{i=1}^{\lfloor m/2\rfloor}\Xi\left(\frac{a_i}{a_{m+1-i}}\right)$$
where $g=kak'$ is the Cartan decomposition of $g$ and $a=\diag{a_1,\dots,a_m}$ is such that $|a_i|\geq|a_j|$ for $i>j$ and $a_i=p^{k_i}$. The function $\Xi$ is the Harish-Chandra function of $\on{PGL}_2(\bQ_p)$, where we write $\Xi(x)=\Xi(\diag{x,1})$ and is explicitly calculated to be $$\Xi(p^k)=\frac{1}{p^{k/2}}\frac{k(p-1)+p+1}{p+1}.$$ The theorem now follows immediately after we sacrifice part of the exponent to get rid of the linear term appearing in $\Xi$.
\ep\medbreak
Let us mention that it is also possible to give an elementary proof similar to \cite{howetan} where the real group $\SLR$ is treated. For the rank $1$ group $\SLQ[2]$ one can use the fact that $G/K$ defines geometrically a $p+1$-regular tree $\cT$. The eigenvalues of the associated Laplacian parametrize the irreducible unitary representations of $G$. In particular, $X/K$ is a finite graph and thus has only finitely many representations appearing. We paraphrase
\begin{theo}
\label{decay2}
There holds an analogous statement for $\SLQ[2]$ as in Theorem~\ref{decay3} but with decay rate now depending on $\Gamma$.
\end{theo}

\acks
We want to thank the anonymous referee for pointing out a simplification that appeared in \cite{shirali} and improve the exponent of Theorem \ref{theorem}. This is part of the author's PhD-thesis supervised by Prof.\ Einsiedler to whom we are grateful not only for proposing this problem to us, but also for providing the dynamical and representational tools we learned from him to study it.


\begin{thebibliography}{CHH88}

\bibitem[AW67]{adler}
R.~L. Adler and B.~Weiss.
\newblock Entropy, a complete metric invariant for automorphisms of the torus.
\newblock {\em Proc. Nat. Acad. Sci. U.S.A.}, 57:1573--1576, 1967.

\bibitem[Bow71]{bowen}
Rufus Bowen.
\newblock Entropy for group endomorphisms and homogeneous spaces.
\newblock {\em Trans. Amer. Math. Soc.}, 153:401--414, 1971.

\bibitem[BW00]{wallach}
A.~Borel and N.~Wallach.
\newblock {\em Continuous cohomology, discrete subgroups, and representations
  of reductive groups}, volume~67 of {\em Mathematical Surveys and Monographs}.
\newblock American Mathematical Society, Providence, RI, second edition, 2000.

\bibitem[CHH88]{chh}
M.~Cowling, U.~Haagerup, and R.~Howe.
\newblock Almost {$L^2$} matrix coefficients.
\newblock {\em J. Reine Angew. Math.}, 387:97--110, 1988.

\bibitem[Cow79]{cowling}
Michael Cowling.
\newblock Sur les coefficients des repr\'esentations unitaires des groupes de
  {L}ie simples.
\newblock In {\em Analyse harmonique sur les groupes de {L}ie ({S}\'em.,
  {N}ancy-{S}trasbourg 1976--1978), {II}}, volume 739 of {\em Lecture Notes in
  Math.}, pages 132--178. Springer, Berlin, 1979.

\bibitem[CT91]{cover}
Thomas~M. Cover and Joy~A. Thomas.
\newblock {\em Elements of information theory}.
\newblock Wiley Series in Telecommunications. John Wiley \& Sons Inc., New
  York, 1991.
\newblock A Wiley-Interscience Publication.

\bibitem[EKL06]{EKL}
Manfred Einsiedler, Anatole Katok, and Elon Lindenstrauss.
\newblock Invariant measures and the set of exceptions to {L}ittlewood's
  conjecture.
\newblock {\em Ann. of Math. (2)}, 164(2):513--560, 2006.

\bibitem[EL10]{pisa}
M.~Einsiedler and E.~Lindenstrauss.
\newblock Diagonal actions on locally homogeneous spaces.
\newblock In {\em Homogeneous flows, moduli spaces and arithmetic}, volume~10
  of {\em Clay Math. Proc.}, pages 155--241. Amer. Math. Soc., Providence, RI,
  2010.

\bibitem[EW11]{vol1}
Manfred Einsiedler and Thomas Ward.
\newblock {\em Ergodic theory with a view towards number theory}, volume 259 of
  {\em Graduate Texts in Mathematics}.
\newblock Springer-Verlag London Ltd., London, 2011.

\bibitem[HM79]{hm}
Roger~E. Howe and Calvin~C. Moore.
\newblock Asymptotic properties of unitary representations.
\newblock {\em J. Funct. Anal.}, 32(1):72--96, 1979.

\bibitem[HT92]{howetan}
Roger Howe and Eng-Chye Tan.
\newblock {\em Nonabelian harmonic analysis}.
\newblock Universitext. Springer-Verlag, New York, 1992.
\newblock Applications of ${{\rm{S}}L}(2,{{\bf{R}}})$.

\bibitem[Kad14]{shirali}
Shirali Kadyrov.
\newblock Effective uniqueness of Parry measure and exceptional sets in ergodic theory.
\newblock {\em Monatshefte f\"ur Mathematik}, to appear (2014).

\bibitem[Mar91]{margulisbook}
G.~A. Margulis.
\newblock {\em Discrete subgroups of semisimple {L}ie groups}, volume~17 of
  {\em Ergebnisse der Mathematik und ihrer Grenzgebiete (3) [Results in
  Mathematics and Related Areas (3)]}.
\newblock Springer-Verlag, Berlin, 1991.

\bibitem[MT94]{margulis}
G.~A. Margulis and G.~M. Tomanov.
\newblock Invariant measures for actions of unipotent groups over local fields
  on homogeneous spaces.
\newblock {\em Invent. Math.}, 116(1-3):347--392, 1994.

\bibitem[Oh02]{oh}
Hee Oh.
\newblock Uniform pointwise bounds for matrix coefficients of unitary
  representations and applications to {K}azhdan constants.
\newblock {\em Duke Math. J.}, 113(1):133--192, 2002.

\bibitem[Pol11]{Polo}
Fabrizio Polo.
\newblock {\em Equidistribution in chaotic dynamical systems}.
\newblock ProQuest LLC, Ann Arbor, MI, 2011.
\newblock Thesis (Ph.D.), The Ohio State University.

\bibitem[PR94]{platonov}
Vladimir Platonov and Andrei Rapinchuk.
\newblock {\em Algebraic groups and number theory}, volume 139 of {\em Pure and
  Applied Mathematics}.
\newblock Academic Press, Inc., Boston, MA, 1994.
\newblock Translated from the 1991 Russian original by Rachel Rowen.

\bibitem[Rat98]{ratner}
Marina Ratner.
\newblock On the {$p$}-adic and {$S$}-arithmetic generalizations of
  {R}aghunathan's conjectures.
\newblock In {\em Lie groups and ergodic theory ({M}umbai, 1996)}, volume~14 of
  {\em Tata Inst. Fund. Res. Stud. Math.}, pages 167--202. Tata Inst. Fund.
  Res., Bombay, 1998.

\bibitem[Ser92]{serre}
Jean-Pierre Serre.
\newblock {\em Lie algebras and {L}ie groups}, volume 1500 of {\em Lecture
  Notes in Mathematics}.
\newblock Springer-Verlag, Berlin, second edition, 1992.
\newblock 1964 lectures given at Harvard University.

\bibitem[Sha00]{shalom}
Yehuda Shalom.
\newblock Rigidity, unitary representations of semisimple groups, and
  fundamental groups of manifolds with rank one transformation group.
\newblock {\em Ann. of Math. (2)}, 152(1):113--182, 2000.

\bibitem[Sil79]{silberger}
Allan~J. Silberger.
\newblock {\em Introduction to harmonic analysis on reductive {$p$}-adic
  groups}, volume~23 of {\em Mathematical Notes}.
\newblock Princeton University Press, Princeton, N.J., 1979.
\newblock Based on lectures by Harish-Chandra at the Institute for Advanced
  Study, 1971--1973.

\end{thebibliography}
\end{document}